\documentclass[french,leqno,a4paper, twoside, 12pt]{article}   
\usepackage[french]{babel} 
\usepackage[T1]{fontenc}
\usepackage[latin1]{inputenc}
\usepackage{amsfonts,amsmath,amssymb,amscd,amsthm}
\usepackage{enumerate,paralist,mathtools}
\usepackage[usenames,dvipsnames]{xcolor}
\usepackage[top=40pt, bottom=40pt, left=40pt, right=40pt]{geometry}

\usepackage{enumerate} 
\usepackage[pdftex]{hyperref} 
\hypersetup{
    bookmarks=true,         
    unicode=false,          
    pdftoolbar=true,        
    pdfmenubar=true,        
    pdffitwindow=false,     
    pdfstartview={FitH},    
    pdftitle={My title},    
    pdfauthor={Author},     
    pdfsubject={Subject},   
    pdfcreator={Creator},   
    pdfproducer={Producer}, 
    pdfkeywords={keyword1} {key2} {key3}, 
    pdfnewwindow=true,      
    colorlinks=true,       
    linkcolor=blue,          
    citecolor=red,        
    filecolor=magenta,      
    urlcolor= black
}

\usepackage{tikz}  
\usetikzlibrary{matrix,arrows}

\usepackage{alltt} 
\usepackage{verbatim}

\usepackage{makeidx} 

\newtheorem{thm}[equation]{Théorème}

\theoremstyle{definition}

\swapnumbers

\newcommand{\N}{\mathbb{N}}
\newcommand{\CC}{\mathbb{C}}

\newcommand{\D}{\mathbb{D}}

\newcommand{\DD}{\partial \mathbb{D}}
\newcommand{\T}{\mathbb{T}}

\makeindex
\begin{document}
\title{Une note sur l'article de S. Takahashi}

\author{TRAN Duc-Anh \footnote{Department of Mathematics, Hanoi National University of
Education, 136 Xuan Thuy St., Hanoi, Vietnam
\newline \href{mailto:ducanh@hnue.edu.vn}{ducanh@hnue.edu.vn}}}\date{}

\maketitle

\abstract Nous démontrons que  la condition nécessaire et suffisante
de S. Takahashi pour le problème d'interpolation de
Nevanlinna-Pick-Carathéodory peut être obtenue de la théorie de D.
Sarason. \endabstract

\section{Enoncé du problème d'interpolation de Nevanlinna-Pick-Carathéodory}
Il faut dire que cette dénomination de ce problème d'interpolation
n'est pas standard et en réalité, dans l'article de S.
Takahashi\cite{Takahashi89}, elle l'appelle par un autre nom. Voici
l'énoncé du problème.

Soit $\D$ le disque unité dans $\CC.$ Soient
$\alpha_1,\alpha_2,\ldots, \alpha_N$ N points distincts dans $\D.$ A
chaque $\alpha_i,$ on associe $n_i$ nombres complexes
$c_{i0},c_{i1},\ldots, c_{i,n_i-1}.$

\paragraph{Problème} \textit{Quelles sont les conditions pour qu'il existe une fonction
holomorphe $\phi\colon \D\to \D$ telle que \begin{equation}\phi(z) =
\sum_{j=0}^{n_i-1}(z-\alpha_i)^j + O((z-\alpha_i)^{n_i})\hskip1cm
(i=1\ldots N)?\label{probleme}\end{equation}}

\section{Condition de S. Takahashi}
Dans cette section, nous citons la condition nécessaire et
suffisante dans l'article de S. Takahashi\cite{Takahashi89}.

D'abord, on associe aux deux points $\alpha_i$ et $\alpha_j$ la
matrice $\Gamma_{ij}$ de façon suivante. On considère le
développment suivant $$\frac{1}{1-z\bar{\zeta}} = \sum_{k,l\geq
0}a_{kl}(z-\alpha_i)^k(\overline{\zeta-\alpha_j})^l.$$ Posons alors
$$\Gamma_{ij} = \begin{pmatrix} a_{00} &\cdots &a_{0,n_j-1}\\
\vdots&&\vdots \\ a_{n_i-1,0} &\ldots& a_{n_i-1,n_j-1}
\end{pmatrix}, \hskip1cm \Gamma =\begin{pmatrix}
  \Gamma_{11}&\cdots&\Gamma_{1N} \\ &\cdots\cdots& \\
  \Gamma_{N1}&\cdots&\Gamma_{NN}
\end{pmatrix}.$$ Ensuite, posons $$C_i = \begin{pmatrix}  c_{i0} & &\\
\vdots&\ddots&\\ c_{i,n_i-1} &\cdots&c_{i0}
\end{pmatrix}, \hskip1cm C=\begin{pmatrix}
  C_1 &&\\ &\ddots& \\ &&C_N
\end{pmatrix},$$ et finalement \begin{equation}\label{condition_T}\boxed{A= \Gamma - C\cdot \Gamma
\cdot C^{\ast}}\end{equation} où $C^{\ast}$ est la matrice
transposée conjuguée de $C.$ Remarquons que $A$ est une matrice
carrée hermitienne de taille $n_1+n_2+\ldots +n_N.$

\begin{thm}[Condition de S. Takahashi]Pour qu'il existe une fonction
$\phi$ qui satisfait à \eqref{probleme}, la condition nécessaire et
suffisante est que la matrice $A$ ci-dessus est positive
semi-définie.
\end{thm}

\section{Une autre démonstration en utilisant la théorie de D. Sarason}

Dans cette section, nous montrons que la condition de S. Takahashi
est une conséquence de la théorie de D. Sarason \cite{Sarason66}.
Nous commençons par quelques notations.

Notons $H^2$ et $H^{\infty}$ les espaces de Hardy sur le disque
$\D.$ Notons $S\colon H^2\to H^2$ le décalage à droite, i.e.
$(Sf)(z)=zf(z)$ pour toute $f\in H^2.$ Pour chaque $\phi \in
H^{\infty},$ notons $M_{\phi}\colon H^2\to H^2$ l'opérateur de
multiplication par $\phi,$ i.e. $M_{\phi}f =\phi f$ pour toute $f\in
H^2,$ et notons $M_{\phi}^{\ast}$ l'adjointe de $M_{\phi}.$ Pour $f,
g\in H^2,$ $\langle f,g\rangle $ désigne le produit scalaire dans
$H^2.$ Et finalement, posons
$$\psi(z) = \prod_{i=1}^N
\left(\frac{z-\alpha_i}{1-\bar{\alpha}_iz}\right)^{n_i}$$ pour $z\in
\D.$ Remarquons que $\psi$ est une fonction intérieure, et elle
jouera le même rôle que la fonction $\psi$ dans \cite{Sarason66}.

D'abord, on cherche des fonctions $h$ qui permettent de reproduire
les dérivées à certain point du disque $\D,$ i.e. pour $f\in H^2,$
$$\langle f, h\rangle = f^{(m)}(\alpha) $$ pour $\alpha\in \D$ et $m\in \N$ fixés.

Pour cela, on utilise l'intégrale de Cauchy. Par un argument de
limite, on peut supposer que les fonctions de $H^2$ sont définies
sur le cercle unité $\mathbb{T}.$ Soient $f\in H^2,$ $\alpha\in \D$
et $m\in\N.$ Par l'intégrale de Cauchy, on a $$f(\alpha) =
\frac{1}{2\pi i}\int_{\DD} \frac{f(z)}{z-\alpha}dz.$$ On a donc
$$f^{(m)}(\alpha) = \frac{m!}{2\pi i}\int_{\DD}\dfrac{f(z)}{(z-\alpha)^{m+1}}dz.$$

On transforme cette intégrale en une intégrale de Lebesgue sur le
cercle $\T$ comme suit \begin{align*}   \frac{f^{(m)}(\alpha)}{m!} &
= \frac{1}{2\pi
i}\int_{0}^{2\pi}\dfrac{f(e^{i\theta})}{(e^{i\theta}-\alpha)^{m+1}}ie^{i\theta}d\theta
\\ & = \frac{1}{2\pi} \int_{0}^{2\pi}\dfrac{f(e^{i\theta})}{(1-e^{-i\theta}\alpha)^{m+1}}e^{-(m+1)i\theta}e^{i\theta}d\theta
\\ & = \frac{1}{2\pi} \int_{0}^{2\pi}\dfrac{f(e^{i\theta})}{(1-e^{-i\theta}\alpha)^{m+1}}e^{-mi\theta}d\theta
\\ & = \left\langle f, \frac{z^m}{(1-\bar{\alpha}z)^{m+1}}\right\rangle.
\end{align*}

Posons alors $$k_{\alpha, m}(z) =
\frac{z^m}{(1-\bar{\alpha}z)^{m+1}}$$ et il est facile de voir que
$k_{\alpha,m}$ permet de reproduire le coefficient du terme
$(z-\alpha)^m$ du développement de Taylor de $f\in H^2.$

Posons ensuite $$K = \mathrm{span}\{k_{\alpha_i, m}~:~ 0\leq m\leq
n_i-1, 1\leq i\leq N\}$$ et il est facile de voir que $\dim K = n_1
+\ldots +n_N$ et  que $K = H^2\ominus \psi H^2.$

Maintenant, supposons que $\phi$ est la fonction interpolant
\eqref{probleme}. Il nous faut calculer l'action de
$M_{\phi}^{\ast}$ sur $K.$ Pour cela, prenons une fonction test
$f\in H^2,$ et on a
\begin{align*}
  \langle f, M^{\ast}_{\phi} k_{\alpha,m}\rangle &= \dfrac{(\phi
  f)^{(m)}(\alpha)}{m!} = \sum_{i=0}^m\dfrac{C^i_m \phi^{(m-i)}(\alpha)
  f^{(i)}(\alpha)}{m!} \\ & = \sum_{i=0}^m
  \dfrac{\phi^{(m-i)}(\alpha)}{(m-i)!} \dfrac{f^{(i)}(\alpha)}{i!}
  \\ & = \left\langle f, \sum_{i=0}^m \overline{\dfrac{\phi^{(m-i)}(\alpha)}{(m-i)!}}\cdot k_{\alpha,
  i}\right\rangle.
\end{align*}

Donc $$ M^{\ast}_{\phi} k_{\alpha,m} = \sum_{i=0}^m
\overline{\dfrac{\phi^{(m-i)}(\alpha)}{(m-i)!}}\cdot k_{\alpha,
  i}.$$

Posons $$D_i = \begin{pmatrix} \bar{c}_{i0} & \bar{c}_{i1}&\cdots
&\bar{c}_{i,n_i-2} & \bar{c}_{i,n_i-1} \\ &\bar{c}_{i0}&
\bar{c}_{i1} &\cdots & \bar{c}_{i,n_i-2} \\ &&\ddots& \ddots
&\vdots\\ &&&\ddots&\vdots \\ &&&&\bar{c}_{i0}
\end{pmatrix}\mbox{ et } D = \begin{pmatrix}
  D_1 & & \\ &\ddots& \\ && D_N
\end{pmatrix}.$$ On peut voir que $D$ est la matrice de
$M_{\phi}^{\ast}\colon K\to K$  dans la base $\{k_{\alpha_i, m}~:~
1\leq i\leq N, 0\leq m\leq n_i-1\}.$

Maintenant, on cherche des conditions pour qu'il existe $\phi$
interpolant \eqref{probleme}. Posons $T\colon K\to K$ un
endomorphisme qui satisfait au fait que $$T k_{\alpha_i, m} =
\sum_{k=0}^m \bar{c}_{ik}\cdot k_{\alpha_i,
  k},$$ autrement dit que $T$ a pour matrice représentante $D$ dans la
  base $\{k_{\alpha_i, m}~:~
1\leq i\leq N, 0\leq m\leq n_i-1\}.$ Il n'est pas difficile de voir
la commutativité de $T$ avec le décalage à gauche $S^{\ast}\colon
K\to K.$

\paragraph{Condition d'interpolation} \textit{ Par la théorie de D. Sarason\cite{Sarason66}, la condition
nécessaire et suffisante pour l'existence d'une fonction $\phi$
interpolant \eqref{probleme} est que $T$ doit être une contraction,
cela dit $1 - T^{\ast}T\geq 0.$}

 C'est équivalent au caractère
positif semi-défini de la matrice suivante
$$\left[\langle (1-T^{\ast}T)v_i, v_j\rangle\right]_{1\leq i, j\leq n_1+\ldots +n_N}
$$  avec $\{v_1,\ldots, v_{n_1+\ldots + n_N}\}$ certaine base de
$K.$ C'est-à-dire que $$\left[\langle v_i,v_j\rangle - \langle
Tv_i,Tv_j\rangle\right]_{1\leq i,j\leq n_1+\ldots +n_N}\geq 0.$$
Soit $G$ la matrice de Gram de la base $\{k_{\alpha_i, m}~:~ 1\leq
i\leq N, 0\leq m\leq n_i-1\}.$

La condition nécessaire et suffisante est donc
\begin{equation}\label{condition_S}\boxed{G - D^T G \bar{D}\geq 0}.\end{equation} On
précise les entrées de $$G =
\begin{pmatrix} G_{11} &\cdots & G_{1N} \\ &\cdots\cdots& \\
G_{N1}&\cdots&G_{NN}
\end{pmatrix}$$ où

 $$G_{ij} = \begin{pmatrix}  \langle k_{\alpha_i, 0},
k_{\alpha_j,0} \rangle &\cdots & \langle k_{\alpha_i, 0},
k_{\alpha_j, n_j-1}\rangle \\ \vdots && \vdots \\ \langle
k_{\alpha_i, n_i-1},k_{\alpha_j,0} \rangle &\cdots& \langle
k_{\alpha_i,n_i-1},k_{\alpha_j,n_j-1}\rangle
\end{pmatrix}$$ pour tous $1\leq i, j\leq N.$

\section{Comparaison avec la condition de S. Takahashi}

On voit immédiatement que $D = C^{\ast}.$ Maintenant, on montre que
$G = \bar{\Gamma}.$

Pour faire cela, il faut calculer le produit scalaire suivant
\begin{align*} \langle k_{\alpha, m}, k_{\beta, n} \rangle & =
\left\langle \frac{z^m}{(1-\bar{\alpha}z)^{m+1}},
\frac{z^n}{(1-\bar{\beta}z)^{n+1}} \right\rangle \\ & = \frac{1}{n!}
\left.\dfrac{d^n}{dz^n}\right|_{z=\beta}
\dfrac{z^m}{(1-\bar{\alpha}z)^{m+1}}.
\end{align*} Quand $\alpha =\alpha_i$ et $\beta=\alpha_j,$ cette
valeur est l'entrée en position $(m,n)$ de $G_{ij}.$

Ensuite, considérons l'entrée de même position de $\Gamma_{ij}.$
C'est le coefficient $a_{mn}$ dans le développement
$$\frac{1}{1-z\bar{\zeta}} =\sum_{k,l\geq 0}a_{kl}(z-\alpha)^k
(\overline{\zeta-\beta})^l.$$ On remarque que
\begin{align*}\sum_{i\geq 0}a_{mi}(\overline{\zeta-\beta})^i & =
\frac{1}{m!} \left.\frac{d^m}{dz^m}\right|_{z=\alpha}
\frac{1}{1-z\bar{\zeta}} \\ & =
\left.\frac{\bar{\zeta}^m}{(1-z\bar{\zeta})^{m+1}}\right|_{z=\alpha}
= \frac{\bar{\zeta}^m}{(1-\alpha \bar{\zeta})^{m+1}}.\end{align*}
Donc $$a_{mn}= \frac{1}{n!}
\left.\frac{d^n}{d\bar{\zeta}^n}\right|_{\zeta=\beta}
\frac{\bar{\zeta}^m}{ (1-\alpha\bar{\zeta})^{m+1}} =
\overline{\langle k_{\alpha, m}, k_{\beta, n} \rangle}.$$ On déduit
que $G=\bar{\Gamma}$ et que $$G- D^T G\bar{D} = \overline{\Gamma -
C\Gamma C^{\ast}},$$ i.e. les deux conditions sont identiques.

\paragraph{Remerciement} Nous remercions Monsieur Pascal Thomas pour
nous envoyer l'article de S. Takahashi.

\end{document}